\documentclass[preprint,12pt]{elsarticle}
\usepackage{amssymb}
\usepackage{amsmath}
\usepackage{graphicx}
\usepackage{epsfig}
\usepackage{color}
\usepackage{amsthm}
\usepackage{verbatim,subfigure,rotating,slashbox,multirow,booktabs}

\usepackage{url}

\newtheorem{varthm}{Theorem}
\newtheorem{theorem}{Theorem}
\newtheorem{corollary}{Corollary}
\newtheorem{lemma}{Lemma}

\newtheorem{proposition}{Proposition}

\newcommand{\lref}[1]{Lemma \ref{lem:#1}}
\newcommand{\cref}[1]{Corollary \ref{cor:#1}}

\newcommand{\sref}[1]{Section \ref{sec:#1}}
\newcommand{\tref}[1]{Theorem \ref{thm:#1}}

\newcommand{\tabref}[1]{Table \ref{tab:#1}}
\newcommand{\fref}[1]{Figure \ref{fig:#1}}

\newcommand{\T}{{\textsf{T}}}
\newcommand{\X}{{\textsf{X}}}
\newcommand{\fx}{x_f}
\newcommand{\fy}{y_f}
\newcommand\set[1]{\{ #1\}}
\newcommand\flr[2]{\left\lfloor \frac{#1}{#2} \right\rfloor}
\newcommand\cln[2]{\left\lceil \frac{#1}{#2} \right\rceil}

\newcommand\ldown[1]{l_{#1}^\downarrow}
\newcommand\lup[1]{l_{#1}^\uparrow}
\newcommand\rdown[1]{r_{#1}^\downarrow}
\newcommand\rup[1]{r_{#1}^\uparrow}

\newcommand\tright[1]{t_{#1}^\rightarrow}
\newcommand\tleft[1]{t_{#1}^\leftarrow}
\newcommand\bright[1]{b_{#1}^\rightarrow}
\newcommand\bleft[1]{b_{#1}^\leftarrow}

\newcommand\tpl{\mathbf{tl}}
\newcommand\tpr{\mathbf{tr}}
\newcommand\btl{\mathbf{bl}}
\newcommand\btr{\mathbf{br}}

\newcommand{\Tn}[1]{\mathbf T_{#1}}

\journal{IWOCA'11 Special Issue}

\begin{document}
\begin{frontmatter}
  \title{Monomer-dimer tatami tilings of square regions}

  \author{Alejandro Erickson}
  \address{Department of Computer Science, University of Victoria, V8W
    3P6, Canada}
  \ead{ate@uvic.ca}
  \ead[url]{http://alejandroerickson.com}
  \author{Mark Schurch}
\address{Mathematics and Statistics, University of Victoria, V8W 3R4, Canada}
\ead{mschurch@uvic.ca}

\begin{abstract}
  We prove that the number of monomer-dimer tilings of an $n\times n$
  square grid, with $m<n$ monomers in which no four tiles meet at any
  point is $m2^m+(m+1)2^{m+1}$, when $m$ and $n$ have the same parity.
  In addition, we present a new proof of the result that there are
  $n2^{n-1}$ such tilings with $n$ monomers, which divides the tilings
  into $n$ classes of size $2^{n-1}$.  The sum of these tilings over
  all monomer counts has the closed form $2^{n-1}(3n-4)+2$ and,
  curiously, this is equal to the sum of the squares of all parts in
  all compositions of $n$. We also describe two algorithms and a Gray
  code ordering for generating the $n2^{n-1}$ tilings with $n$
  monomers, which are both based on our new proof.
\end{abstract}

\end{frontmatter}

\section{Introduction}
\label{sec:intro}

Tatami mats are a traditional Japanese floor covering, whose outside
is made of soft woven rush straw, and whose core is stuffed with rice
straw.  They have been in use by Japanese aristocracy since the $12$th
century, and are so integral to Japanese culture, that they are often
used to describe the square footage of a room.  A standard full mat
measures $6'\times 3'$, and a half mat is $3'\times 3'$.

An arrangement of the mats in which no four meet at any point is often
preferred.  We call such arrangements monomer-dimer \emph{tatami
  tilings}.  The monomer-dimer tiling in \fref{nonTatami} violates the
tatami condition, and the one in \fref{example} does not.
\begin{figure}[h]
  \centering
  \includegraphics{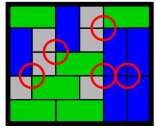}
  \caption{The tatami condition is that no four tiles may meet at a
    point.  All such
violations are circled.}
\label{fig:nonTatami}
\end{figure}

\begin{figure}[h]
  \begin{center}
    \includegraphics[height = 1in]{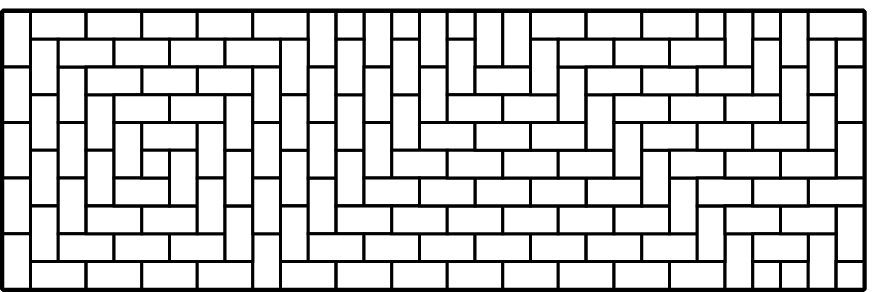}      
  \end{center}
  \caption{ A tatami tiling of the $10\times 31$ grid with $10$
    monomers, $64$ vertical dimers and $86$ horizontal dimers.
    Compare this with \fref{comprehensiveTiling}.}
  \label{fig:example}
\end{figure}

Tilings with polyominoes are well studied, and they appear in physical
models, the theory of regular languages, and combinatorics (see
\cite{BenedettoLoehr2008,Stanley1985,GaleGolombHaas1996,MerliniSprugnoliVerri2000,HockMcQuistan1984}).

It is natural to ask, given a room's size, how many different
arrangements of mats are there?  Motivated by an exercise posed by
Knuth (\cite{Knuth2011}), Ruskey and Woodcock published the ordinary
generating functions for fixed height, dimer-only tatami tilings in
\cite{RuskeyWoodcock2009}.

In Fall 2009 our research group discovered the structure which is
published in \cite{EricksonRuskeySchurch2011}, and described in
\sref{structure}.  Problems such as the enumeration of tatami tilings
with $r$ rows, $c$ columns, and $m$ monomers, have since become within
reach.  This paper completes such an enumeration for square grids.






Alhazov et al. followed up on the aforementioned dimer-only research
with a treatment of odd-area tatami tilings which include a single
monomer.  They closed \cite{ALHAZOVMORITAIWAMOTO2010} with this remark:
\begin{quote}
  However, the variety of tilings with arbitrary number of
  monominoes is quite ``wild'' in sense that such tilings cannot be
  easily decomposed, see Figure 11; therefore, most results presented
  here do not generalize to arbitrary number of monominoes, the
  techniques used here are not applicable, and it is expected that any
  characterization or enumeration of them would be much more
  complicated.
\end{quote}

The structure we found, however, reveals the opposite; the tilings
with an arbitrary number of monomers \emph{are} easily decomposed.
The decomposition has a satisfying symmetry, it is ammenable to
inductive arguments, and it shows that the complexity of a tatami
tiling is linear in the dimensions of the grid (compare \fref{example}
with \fref{comprehensiveTiling}).  We use it extensively to prove our
main result.  Let $T(n,m)$ be the number of $n\times n$ tatami tilings
with exactly $m$ monomers.
\begin{varthm}
\label{varthm:square}
  If $n$ and $m$ have the same parity, and $m<n$, then $T(n,m) = m2^{m} +
  (m+1)2^{m+1}$.
\end{varthm}
This confirms Conjecture 3 in \cite{EricksonRuskeySchurch2011} for $d=0$ (as was promised
there), and completes the enumeration of tatami tilings of square
grids, since $T(n,m)=0$ when $m>n$.

The proof that $T(n,n)=n2^{n-1}$, in \cite{EricksonRuskeySchurch2011},
counts tilings recursively. Due to the symmetry of the solution it is
natural to search for a way to divide the tilings into $n$ classes of
size $2^{n-1}$; indeed, Knuth mentioned this in a correspondence with
Ruskey (\cite{EricksonSchurch2011}). We remedy the lack of symmetry in
the original proof by presenting a new one which includes the desired
partition and counts the tilings directly, rather than recursively.

Nice round numbers such as these tend to be easily identified in
output generated by exhaustively listing all tilings on a given size
of grid.  In a personal communication Ruskey
(\cite{EricksonSchurch2011}) Knuth noted there appeared to be
$2^{n-1}(3n-4)+2$ tatami tilings of the $n\times n$ grid. We prove
this result, and show that it is equal to the sum of the squares of
all parts in all compositions of $n$.

In \sref{combalg} we use the partition in \tref{tnnn} to exhaustively
generate all $n\times n$ tilings with $n$ monomers, for even $n$, in a
constant amount of time per tiling.  The problem of generating the
tilings in each of these partitions reduces to generating all bit
strings of a fixed length, so our first algorithm is a straightforward
application of this.  The second algorithm generates the tilings in a
Gray code order, where the Gray code operation is the diagonal flip
(defined in \sref{structure}).  This is acheived by generating the
partitions in a certain order, and using a binary reflected Gray code
to generate the bit strings.

\subsection{Structure}
\label{sec:structure}
\begin{figure}[ht]
  \centering%
  \subfigure[A \emph{loner} feature.]{\includegraphics{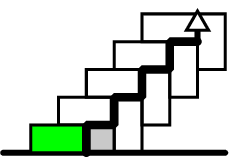}%
    \label{fig:lonersource}}%
  \hspace{0.5in}%
  \subfigure[A \emph{vee} feature.]{\includegraphics{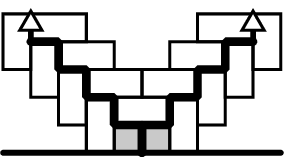}%
    \label{fig:veesource}}%
\caption{These two types of features must have their coloured tiles on
  a boundary, as shown, up to rotation and reflection.}
  \label{fig:lonerveesources}
\end{figure}

\begin{figure}[ht]
  \centering%
  \subfigure[A vertical \emph{bidimer}
  feature.]{\includegraphics{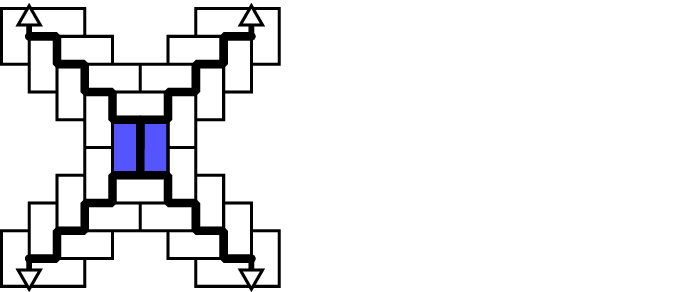}
    \label{fig:bidimersource} }%
  \hspace{0.5in}%
  \subfigure[A
  counterclockwise \emph{vortex}
  feature.]{ \includegraphics{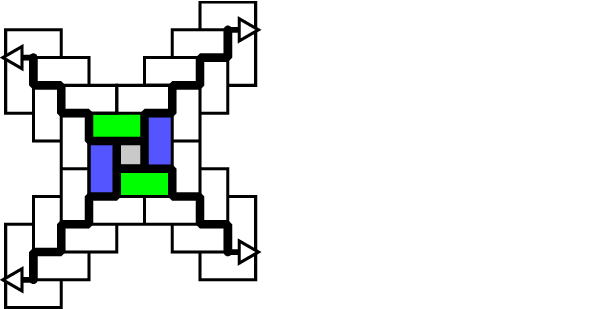} \label{fig:vortexsource}}%
  \caption{These two types of features may appear anywhere in a tiling
    provided that the coloured tiles are within the boundaries of the
    grid.}
  \label{fig:bidimervortexsources}
\end{figure}

Tatami tilings of rectangular grids have an underlying structure,
called the \emph{\T-diagram}, which is an arrangement of four possible
types of \emph{features}, shown in Figures
\ref{fig:lonerveesources}-\ref{fig:bidimervortexsources}, up to
rotation and reflection (see \cite{EricksonRuskeySchurch2011} for a complete explanation).
Each feature consists of a \emph{source}, shown in colour, which
forces the placement of up to four \emph{rays}. Rays propagate to the
boundary of the grid and do not intersect each other.

A valid placement in the grid of a non-zero number of features
determines a tiling uniquely.  Otherwise, the trivial \T-diagram
corresponds to the four possible running bond patterns (brick laying
pattern).

Monomers occur on the boundary of the grid, and inside vortices, but
nowhere else.  The \T-diagram is also a partition of the tiles into
blocks of horizontal and vertical running bond. The tiling in
\fref{comprehensiveTiling} uses all four types of features, displays
this partition into horizontal and vertical running bond. Its
\T-diagram is shown on its own in \fref{tdiagram}.

\begin{figure}[ht]
\centering
 \includegraphics{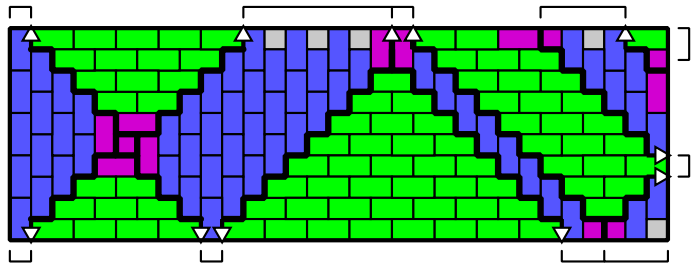}
 \caption{A tiling showing all four types of sources.  Coloured in
   magenta, from left to right they are, a clockwise vortex, a vertical
   bidimer, a loner, a vee, and another loner.}
  \label{fig:comprehensiveTiling}
\end{figure}

\begin{figure}[h]
  \centering
  \includegraphics[width=3in]{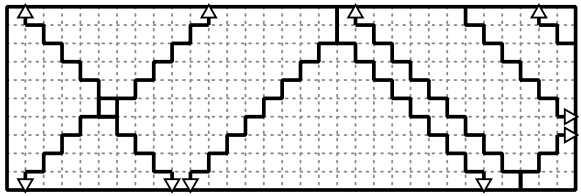}
  \caption{A \T-diagram}
  \label{fig:tdiagram}
\end{figure}

Throughout this paper we consider tilings on the $n\times n$ integer
grid, with the origin at the bottom left of the tiling.  Let the
coordinate of a grid square be the point at its bottom left corner as
well.

A \emph{diagonal} is a rotation of the following: a monomer at
$(x_0,0)$, and (vertical) dimers covering the pairs of grid squares
$(x_0+1+k,k)$ and $(x_0+1+k,k+1)$, for each non-negative $k$ such that
$x_0+1+k\le n-1$.  A diagonal can be \emph{flipped} in place, so that
the monomer moves to $(n-1,n-x_0-1)$, and the dimers change
orientation.  This operation preserves the tatami condition.
Diagonals are used and illustrated frequently in the next section.

\section{Square grids}
\label{sec:square}

\begin{figure}[ht]
  \centering
  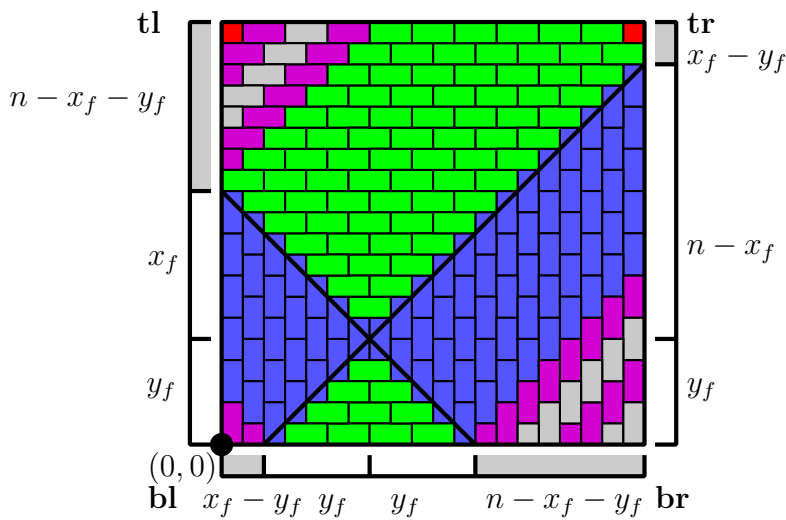
  \caption{Vertical bidimer.}
  \label{fig:vbidimer}
\end{figure}

\begin{figure}[ht]
  \centering
  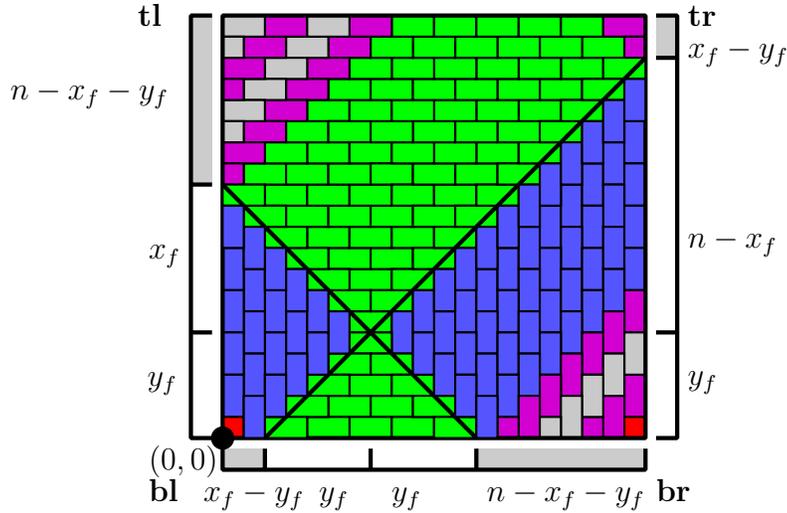
  \caption{Horizontal bidimer.}
  \label{fig:hbidimer}
\end{figure}

\begin{figure}[ht]
  \centering
  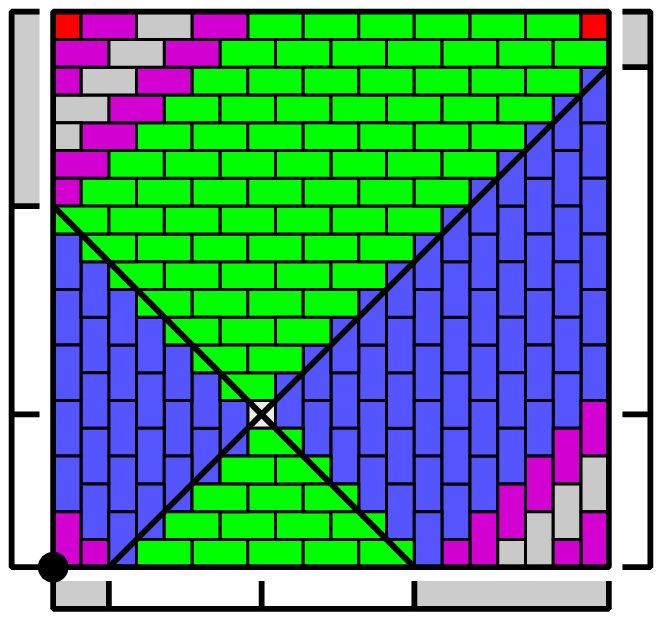
  \caption{Counter clockwise vortex. Note that $\fx$ and $\fy$ are
    not integers.}
    \label{fig:ccvortex}
\end{figure}

Let $\Tn{n,m}$ be the set of $n\times n$ tatami tilings with $m$
monomers.  To find $|\Tn{n,m}|$, we use a setup similar to that of
Theorem 2 in \cite{EricksonRuskeySchurch2011}.  These lemmas on the
composition of such tilings seem apparent from Figures
\ref{fig:vbidimer}--\ref{fig:ccvortex}, but we prove them in general.

We show that any tiling in $\Tn{n,m}$ with $m<n$ has exactly one
bidimer or vortex and we give a one-to-one correspondence between $m$
and the shortest distance from this source to the boundary.  Such a
feature determines all tiles in the tiling except a number of
diagonals that can be flipped independently.  Proving the result
becomes a matter of counting the number of allowable positions for the
bidimer or vortex, each of which contributes a power of $2$ to the
sum.

For example, the $20\times 20$ tiling in \fref{vbidimer} has a
vertical bidimer which forces the placement of the green and blue
tiles, while the remaining diagonals are coloured in alternating grey
and magenta.  There are eight such diagonals, so there are $2^8$
tilings of the $20\times 20$ grid, with a vertical bidimer in the
position shown.  Each of these $2^8$ tilings has exactly $10$
monomers.

By Corollary 1 in \cite{EricksonRuskeySchurch2011}, a tiling in
$\Tn{n,m}$, with $m<n$, must have at least one bidimer or vortex,
which we will call $f$.  Let $(\fx,\fy)$ be the center of $f$ in the
cartesian plane, where the bottom left of the $n\times n$ tiling is at
the origin.

Define $X_f$ as the lines $(x-\fx)+(y-\fy) = 0$ and
$(x-\fx)-(y-\fy)=0$, bounded by the grid boundary, which form an
\X~through $(\fx,\fy)$ (see Figures
\ref{fig:vbidimer}-\ref{fig:ccvortex}).

Without loss of generality, we may re-orient a tiling, by rotating and
reflecting, so that $\fy\le \fx \le n/2$.  The upper arms of $X_f$
intersect the left and right boundaries of the grid, while the lower
arms intersect the bottom.  In this range, the \emph{distance} from
this feature to the boundary of the grid is $\fy$.

\begin{lemma}
  \label{lem:onlyOneX}
  Let $T\in \Tn{n,m}$, with $m<n$.  Then $T$ has exactly one bidimer
  or vortex, $f$, but not both, and no vees.
\end{lemma}
\proof It is easy to see that if $a$ is another bidimer or vortex in
the tiling, then $X_a$ intersects $X_f$, which contradicts the fact
that rays do not cross.  A vee has the same rays as a bidimer, by
replacing the adjacent monomers with a dimer, so the only features
that may appear, besides $f$, are loners.\qed

Let $T_f$ be the tiling with $f$ as its only feature.  This
corresponds to the featurless ``trivial'' tilings in
\cite{EricksonRuskeySchurch2011}.

\begin{lemma}
  \label{lem:ejcLemma3}
  The tiling $T_f$ can be obtained from any tiling containing the
  feature $f$ via a finite sequence of diagonal flips in which each
  monomer moves at most once.  Reversing this sequence gives the
  original tiling.
\end{lemma}
\proof This is a simple modification of the proof of Lemma 3 in
\cite{EricksonRuskeySchurch2011}.\qed

\begin{lemma}
  \label{lem:indepFlips}
  All the diagonals in $T_f$ can be flipped independently.
\end{lemma}
\proof We show that no two diagonals intersect, by sharing a
monomer, or other grid squares.

Let $\alpha$ and $\beta$ be the monomers in distinct diagonals,
possibly with $\alpha = \beta$, which intersect.  Let $L_\alpha$ be
the line through $\alpha$, and let $L_\beta$ be the line through
$\beta$.  Their slopes must differ, because in sharing grid squares,
they must share whole tiles.  Thus we may assume that $L_\alpha$ has
slope $1$, and $L_\beta$ has slope $-1$.  Then $L_\alpha$ and
$L_\beta$ intersect inside the grid, and therefore at least one of
them instersects with $X_f$, which makes the corresponding diagonal
unflippable.  Therefore distinct diagonals may not intersect.\qed

\begin{figure}[ht]
  \centering \hspace{-1.8cm} \subfigure[The range $\fy\le \fx \le n/2$
  is green.]{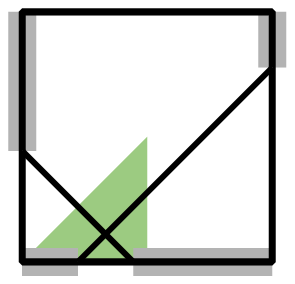 \label{fig:Xa}}%
  \hspace{1cm}
  \subfigure[]{\includegraphics{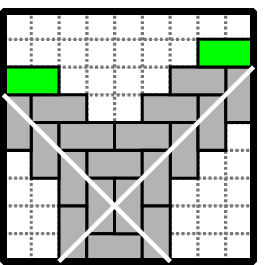} \label{fig:Xvb}}\\
  \subfigure[]{\includegraphics{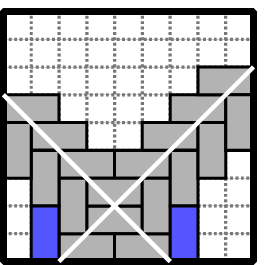} \label{fig:Xhb}}%
  \hspace{0.2cm}
  \subfigure[]{\includegraphics{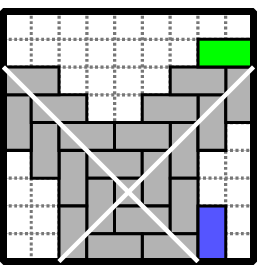} \label{fig:Xcv}}%
  \hspace{0.2cm}
  \subfigure[]{\includegraphics{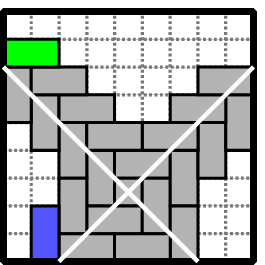} \label{fig:Xccv}}%
  \caption{The segments $\tpl,\tpr,\btl,$ and $\btr$.  When $\tpl$ and
    $\btr$ have positive length, one of their combined tiles is part
    of a ray of $f$, shown in green or blue.  Similarly for $\tpr$
    and $\btl$ combined.}
\end{figure}

The crux of the argument in \tref{square} is this:

\begin{lemma}
  \label{lem:distm}
  If $T_f\in \Tn{n,m}$, with $\fy\le \fx \le n/2$, then $m=n-2\fy$ if
  $f$ is a bidimer, and $m=n-\fy+1$ if $f$ is a vortex.  The number of
  (flippable) diagonals in $T_f$ is

  \begin{align*}
    \left\{\begin{aligned}
        n-2\fy-2&,& \text{ if  } \fy<\fx\text{;}\\
        n-2\fy-1&,& \text{ if  } \fy=\fx<n/2\text{; and}\\
        0&,& \text{ otherwise.}
    \end{aligned}\right.
  \end{align*}
\end{lemma}
\proof Let $\tpl,\tpr,\btl,$ and $\btr$ be the segments of the
boundary of the grid delimited by $X_f$, as shown in grey in
\fref{Xa}.  We count the number of monomers and diagonals in $T_f$ by
measuring the lengths of $\tpl,\tpr,\btl,$ and $\btr$, minus the tiles
in the rays of $f$, casified in Figures \ref{fig:Xvb}--\ref{fig:Xccv}.

Figures \ref{fig:Xvb}--\ref{fig:Xccv} show (in general) that if $\tpl$
and $\btr$ are non-zero, then together they contain exactly one grid
square covered by the ray of $f$.  The same is true of $\tpr$ and
$\btl$.

Because of the above, non-zero $\tpl$ and $\btr$ differ in parity, so
their sum is odd, and therefore they must contain exactly one corner
monomer (which is not in a diagonal).  Similarly for $\tpr$ and
$\btl$.

We tabulate the numbers of monomers and diagonals in $\tpl,\tpr,\btl$
and $\btr$ for different $(\fx,\fy)$ in \tabref{posMons}, and add
pairs of rows to prove each case of the theorem statement.

\begin{table}[ht]
  \centering
\begin{tabular}[h]{llll}
  $(\fx,\fy)$ & Positions & Monomers & Diagonals \\
  \hline
  $\fy =\fx$ & $\btl$ and $\tpr$ & $0$& $0$ \\
  $\fy < \fx$ & $\btl$ and $\tpr$& $\fx-\fy$&$\fx-\fy-1$\\
  $\fy = \fx = n/2$ & $\tpl$ and $\btr$ & $0$ & $0$\\
  $\fy + \fx < n$ & $\tpl$ and $\btr$ & $n-\fx-\fy$&$n-\fx-\fy-1$\\
\end{tabular}
  \caption{}
  \label{tab:posMons}
\end{table}

If $\fy < \fx$, then $m=\fx-\fy + n-\fx-\fy = n-2\fy$ if $f$ is
  a bidimer, and $m=n-2\fy+1$ if $f$ is a vortex.  The number of
  diagonals is $n-2\fy-2$, as required for this case

  If $\fy = \fx < n/2$, then $m=0 + n-\fy - \fy = n-2\fy$ if $f$ is a
  bidimer, and $n-2\fy+1$ if $f$ is a vortex.  The number of diagonals
  is $n-2\fy-1$.

  If $\fy = \fx = n/2$, then $m = 0$ if $f$ is a bidimer, and
  $m=1$ if it is a vortex.  The number of diagonals is $0$.
\qed

We sum up the number of positions which give a particular $m$ in
\lref{distm}, to count all such tilings.

\begin{theorem}\label{thm:square}
  If $n$ and $m$ have the same parity, and $m<n$, then $T(n,m) = m2^{m} +
  (m+1)2^{m+1}$.
\end{theorem}
\proof We count the number of tilings $T_f$ in $\Tn{n,m}$, for each
bidimer or vortex, $f$, satisfying the conditions of \lref{distm}, and
for each of these we count the tilings obtainable via a set of
(independent) diagonal flips.

If $f$ is a bidimer, then $m=n-2k$, where $k$ is the shortest
distance from $(\fx,\fy)$ to the boundary of the grid.  These are the
positions
\begin{align*}
 & (k,k),(k,k+1),(k,k+2),\ldots (k,n-k), (k+1,n-k), \ldots,\\
 & (n-k,n-k),(n-k,n-k-1), \ldots (n-k,k), (n-k-1, k),\ldots (k+1,k),
\end{align*}
of which there are $4(n-k-k)$.  We apply \lref{distm} by making the
necessary rotations and reflections, so that when $m>0$, the four
positions $(k,k), (k,n-k),(n-k,n-k)$ and $(n-k,k)$, have $m-1$
diagonals, while the remaining $4m-4$ positions have $m-2$ diagonals.

The same logic applies when $f$ is a vortex and $m>1$, and the results
are in \tabref{Tnnm}.

By Lemmas \ref{lem:ejcLemma3} and \ref{lem:indepFlips}, all of the
tilings which contain $f$ can be obtained from $T_f$ by making a set
of independent flips.  Thus, the number of these is $2^{d(f)}$, where $d(f)$
is the number of diagonals in $T_f$.

\begin{table}[ht]
  \centering
  \begin{tabular}{llll}
    Type &  Feature  &  Positions &
    Diagonals\\
    \hline
    $\fx=\fy$   & h and v bidimers & $4$ & $m-1$\\
    $\fx=\fy$ & cc and c vortices & $4$ & $m-2$\\
    $\fx < \fy$ & h and v bidimers & $4(m-1)$ & $m-2$ \\
    $\fx<\fy$ & cc and c vortices & $4(m-2)$ & $m-3$
  \end{tabular}
  \caption{Horizontal, vertical, counterclockwise and clockwise
    are abbreviated as h, v, cc
    and c, respectively. We assume $m>0$ if $f$ is a bidimer, and
    $m>1$ if $f$ is a vortex.}
  \label{tab:Tnnm}
\end{table}

\lref{onlyOneX} tells us that there is no other way to obtain a tiling
in $\Tn{n,m}$, so we conclude by summing the $2^{d(f)}$s for each $f$
such that $T_f\in \Tn{n,m}$.  Each term in the following sum comes
from a row of \tabref{Tnnm}, in the same respective order, and
similarly, the three factors in each sum term come from the last three
columns of the table.
\begin{align*}
  T(n,m) =& 2\cdot 4  \cdot 2^{m-1} + 2\cdot 4 \cdot
  2^{m-2} + 2\cdot 4  (m-1) \cdot 2^{m-2} + 2\cdot 4 (m-2)
  \cdot 2^{m-3}\\
  =& 2\cdot 2^{m+1} + 2\cdot 2^m + (m-1)2^{m+1} + (m-2)2^m\\
  =& m2^m + (m+1)2^{m+1}.
\end{align*}
That is, there are $m2^m$ tilings with vortices and $(m+1)2^{m+1}$
with bidimers.  This also satisfies $T(n,0)$ when $f$ is a bidimer, and
$T(n,1)$ when $f$ is a vortex.\qed

This completes the enumeration of $n\times n$ tatami tilings with $m$
monomers. 

We present an alternate proof of Theorem 2 in
\cite{EricksonRuskeySchurch2011}.  The new proof organizes the
$n\times n$ tatami tilings with $n$ monomers into $n$ classes of size
$2^{n-1}$, answering a challenge of Knuth, as noted in
\cite{EricksonSchurch2011}.  In addition to showing a satisfying
symmtery, the resulting partition is used in \sref{combalg} for
exhaustively generating all such tilings in constant amortized time.

The proof in \cite{EricksonRuskeySchurch2011} uses induction on the
dimensions of the grid to show that $n2^{n-1} = 4S(n)$, where
  \begin{align*}
    S(n) = 2^{n-2} + 4 S(n-2), \text{ where } S(1) = \frac{1}{4} \text{
      and } S(2) = 1.
\end{align*}
The recurrence arises from relating the conflicts in flipping pairs of
diagonals in an $(n-2)\times (n-2)$ tiling with those of an $n\times
n$ tiling.

We reuse the setup in \cite{EricksonRuskeySchurch2011}.  Every tatami
tiling of the $n\times n$ grid with $n$ monomers can be obtained by
performing a sequence of flips (of diagonals) on a running bond, in
which no monomer is moved more than once, and the corner monomers
remain fixed.  As an immediate consequence, every tiling has exactly
two corner monomers which are in adjacent corners, and as such, it has
four distinct rotational symmetries.

Let $\Tn{n}$ be the $n\times n$ tilings with $n$ monomers, whose upper
corners have monomers.  It is sufficient to divide the $n2^{n-3}$
tilings of $\Tn{n}$, into $n$ classes of size $2^{n-3}$, and then use
the rotational symmetries of these to obtain the $n$ classes of size
$2^{n-1}$ of all the tilings, $\Tn{n,n}$.

In this context, a \emph{flipped diagonal} or \emph{monomer} refers to one which was
flipped from its initial position in the running bond, in the sequence
of flips referred to above.  Each monomer is in exactly two diagonals
in the trivial tiling, and the monomer is flipped at most once, so its
status can be described by a ternary symbol (see \fref{oddbond}).

\begin{figure}[h]
  \centering
  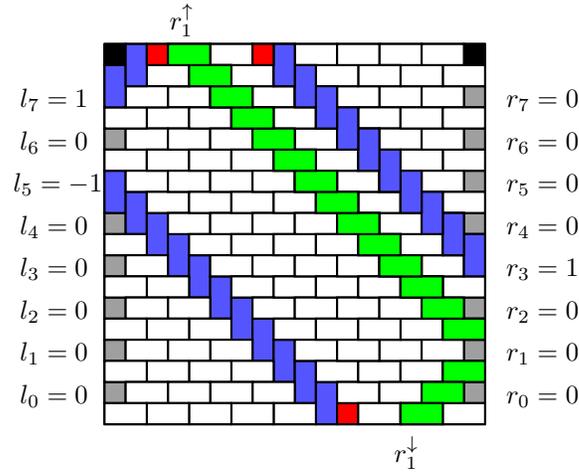
  \caption{This $18\times 18$ tiling is represented by
    $(0,0,0,0,0,-1,0,1)\cdot(0,0,0,1,0,0,0,0)$.  The position of a
    monomer in the trivial tiling reveals the number of dimers in its
    diagonals.  The flipped diagonals are blue, and the diagonals
    $\rup{1}$ and $\rdown{1}$ are green, but they are unflipped.  As
    expected, $d_{18}(\rdown{1}) + d_{18}(\rup{1}) = 3 + 14 = 18-1$.}
  \label{fig:tnnnRep}
\end{figure}

\begin{figure}[h]
  \centering
  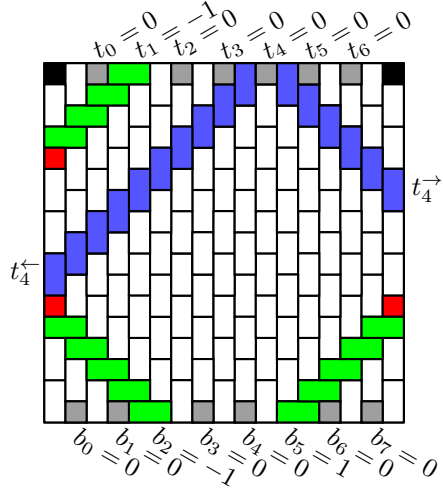
  \caption{This $17\times 17$ tiling is represented by
    $(0,-1,0,0,0,0,0)\cdot(0,0,-1,0,0,1,0,0)$.  The flipped diagonals
    are green, and the two blue diagonals are coloured to show
    $\tleft{4}$ and $\tright{4}$, which are the left and right
    diagonals of the same monomer.  As expected, $d_{17}(\tleft{4}) +
    d_{17}(\tright{4}) = 10 + 6 = 17-1$.}
  \label{fig:tnnnOddRep}
\end{figure}

Putting these ternary symbols in a list provides a representation for
elements of $\Tn{n}$, which contains the status of each monomer as
either unflipped, or flipped in a given direction.  This gives the
sequence of flips required to determine the tiling (See Figures
\ref{fig:tnnnRep}--\ref{fig:tnnnOddRep}).

For even $n$, let $L = (l_0,l_1,\ldots , l_{\frac{n-2}{2}-1})$ be a
ternary string, with $l_i$ representing the $i$th monomer from the
bottom, on the left side of the trivial tiling.  Assign $l_i=0$ if the
monomer is not flipped, $l_i=-1$ if it is flipped down, and $l_i=1$ if
it is flipped up.  The downward diagonal of a monomer is denoted
$\ldown{i}$, and the upward diagonal is $\lup{i}$, regardless of
whether they can, or have been flipped.  Similarly, we define $R =
(r_0,r_1,\ldots, r_{\frac{n-2}{2}-1})$.  Their concatenation, $L\cdot
R$ determines the tiling by describing the diagonal flips, up the left
side, then up the right side.

If $n$ is odd the representation is similar, noting that there are
$\flr{n-2}{2}$ non-fixed monomers along the top and $\cln{n-2}{2}$
along the bottom.  We label from left to right, by $To = (t_0,t_1,
\ldots, t_{\flr{n-2}{2}})$ and $Bo = (b_0,b_1,\ldots,
b_{\cln{n-2}{2}})$, respectively, and let $t_i = -1$ if it is flipped
on its left diagonal.  Similarly, we define $\tright{i}$ as the
rightward diagonal of $t_i$, and so forth.

Let $d_n(a)$ be the number of dimers in the diagonal $a$, also called
the length or size of the diagonal.  It is a function of the index and
direction of $a$:

\begin{align*}
  d_n(a)=\left\{\begin{aligned} 2i+1&,& \text{if }& a=\ldown{i} \text{
        or } a=\rdown{i} \text{
        or } a=\bleft{i};\\
      n-2i-2&,& \text{if }& a=\lup{i} \text{ or } a = \rup{i} \text{
        or } a=\bright{i};\\
      2i+2&,& \text{if }& a=\tleft{i}; \text{ and}\\
      n-2i-3&,& \text{ if }& a= \tright{i}.
    \end{aligned}\right.
\end{align*}

Conflicting pairs of diagonal are identified in a constant number of
operations by using the cases listed in \tabref{conflicts}.  Conflicts
occur only in cases 1(c) and 2(b), and case 1(c) is always a conflict.
In the case of 2(b), there is a conflict if and only if the sizes of
the diagonals sum to at least $n$.

\begin{table}[ht]
  \centering
  \begin{tabular}[h]{ll}
\begin{minipage}[h]{0.5\linewidth}
  A pair of monomers may be
  \begin{enumerate}
\item  on the same side and flipped
  \begin{enumerate}
  \item in the same direction;
  \item away from each other; or
  \item toward each other; or
  \end{enumerate}
\item on different sides and flipped
  \begin{enumerate}
  \item in different directions; or
  \item in the same direction.
  \end{enumerate}
\end{enumerate}
\end{minipage} &
\hspace{0.3cm} \begin{minipage}[h]{0.5\linewidth}
  \begin{tabular}[h]{ll}
    Pair & Conflict 2(b) iff\\
    \hline
    $\ldown{i},\rdown{j}$ &$n \le 2(i+j)+2$\\ 

    $\lup{i}, \rup{j}$& $n\ge 2(i+j)+4$\\ 

    $\bleft{i}, \tleft{j}$& $  n \le 2(i+j) + 3$\\ 

    $\bright{i},\tright{j}$&$n\ge 2(i+j)+5$\\ 
  \end{tabular}
\end{minipage}
  \end{tabular}
  \caption{Diagonal conflict cases.}
  \label{tab:conflicts}
\end{table}

Let $\Tn{n}(a)\subseteq \Tn{n}$, where $a$ is a diagonal (one of
$\ldown{i}, \lup{i}, \rdown{i}, \rup{i}, \bright{i}, \bleft{i},
\tright{i}$, or $\tleft{i}$), be defined as the collection of tilings
in $\Tn{n}$ in which $a$ is the \emph{longest flipped diagonal}; for
each flipped diagonal $b$, distinct from $a$, we have $d_n(b)<d_n(a)$.

Let $\Tn{n}(\varnothing)$ be the set of tilings in which no monomer is
flipped on its longest diagonal.  Note the distinction between a
monomer flipped on its longest diagonal, and the longest flipped
diagonal in the tiling (see \fref{oddbond}).

\begin{figure}[ht]
  \centering
 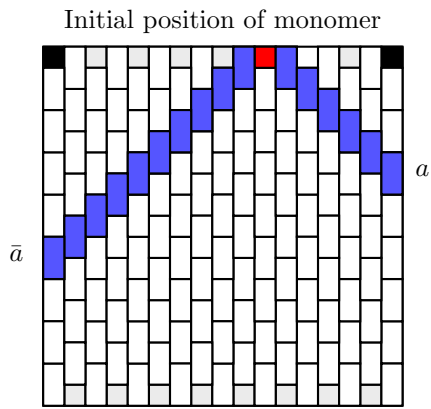
 \caption{Each monomer in the tivial tiling is in two diagonals, $a$
   and $\bar a$.  Note that $d_n(a) + d_n(\bar a) = n-1$}
\label{fig:oddbond}
\end{figure}

Here is the new proof, using the notation defined above.

\begin{theorem}[ \cite{EricksonRuskeySchurch2011}]
  \label{thm:tnnn}
  $T(n,n) = n2^{n-1}$
\end{theorem}
\proof[($n$ classes of size $2^{n-1}$)] Let $D$ be the set of
diagonals in the trivial tiling, and for each monomer in a diagonal
$a$, we denote its other diagonal by $\bar a$ (see \fref{oddbond}).
Let $A = \set{a\in D: d_n(a) > d_n(\bar a)}$ and let $\dot\cup$ denote
a disjoint union.  We prove that
\begin{align*}
\Tn{n} =   \Tn{n}(\varnothing) \dot\cup \left(\dot\bigcup_{a\in A} \Tn{n}(a) \right),
\end{align*}
and, for $a\in A$,
\begin{align}
  \label{eq:setsizes}
  (|\Tn{n}(a)|, |\Tn{n}(\varnothing)|,|A|) = \left\{
    \begin{aligned}
      (2^{n-3},2^{n-2},n-2)&,& \text{ if } n \text{ is even; and}\\
      (2^{n-3}, 2^{n-3},n-1)&,& \text{ if } n \text{ is odd,}
    \end{aligned}\right.
\end{align}
which gives the correct sizes, since $|\Tn{n}| = (n-2)2^{n-3}+2^{n-2}
= (n-1)2^{n-3} + 2^{n-3}$.

An \emph{available diagonal} in $\Tn{n}(a)$ is a diagonal that may be
flipped in a tiling of $\Tn{n}(a)$, subject to the constraints that
$a\in A\cup \set{\varnothing}$ and if $a\neq \varnothing$, then
$d_n(b) < d_n(a)$, and $b$ does not conflict with $a$.

\begin{figure}[h]
  \centering
  \subfigure[$\Tn{17}(\tleft{5})$.]{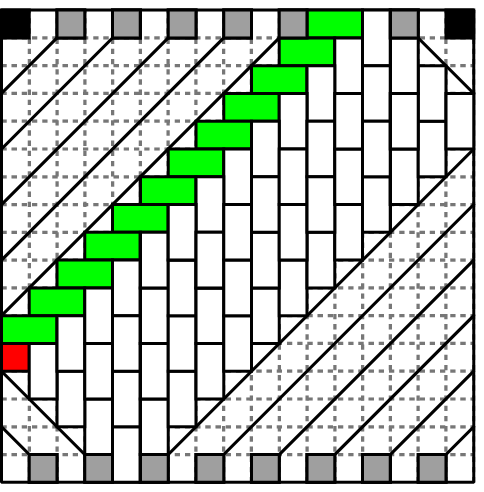 \label{fig:oddlongest}}%
  \hspace{.3cm}%
  \subfigure[ $\Tn{17}(\varnothing)$. Note that $\tleft{3}$ and
  $\tright{3}$ are not available.
  ]{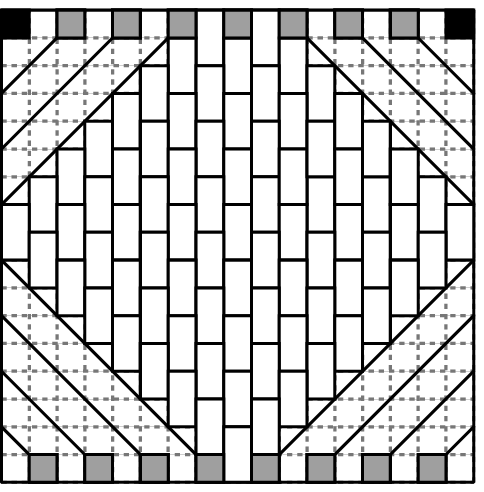 \label{fig:oddnolongest}}%
  \caption{Odd $n$.  Available diagonals are marked by solid lines.}
  \label{fig:oddTnnn}
\end{figure}

\begin{figure}[h]
  \centering
  \subfigure[$\Tn{18}(\ldown{5})$.]{ 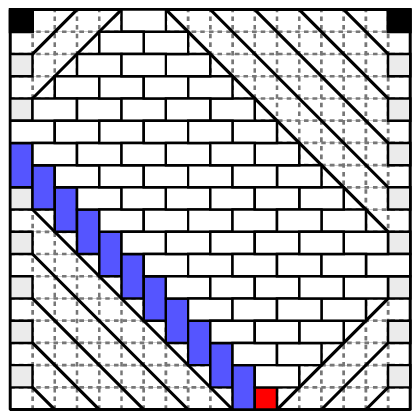 \label{fig:evenlongest}}%
  \hspace{0.3cm}%
  \subfigure[$\Tn{18}(\varnothing)$.]{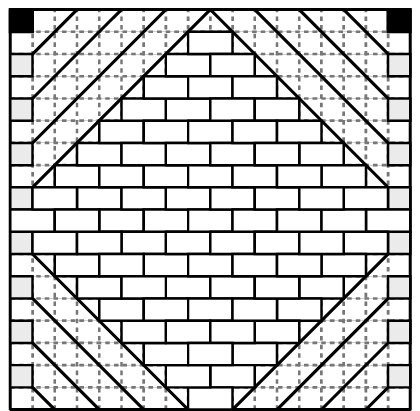 \label{fig:evennolongest}}%
  \caption{Even $n$.  Available diagonals are marked by solid lines.}
  \label{fig:evenTnnn}
\end{figure}

It suffices to show that the available diagonals are pairwise
non-conflicting, and independent, and that for each $\Tn{n}(a)$, there
are $\log_2(|\Tn{n}(a)|)$ of them.  See \tabref{tnnnCases} for all the
cases, and Figures \ref{fig:oddTnnn}--\ref{fig:evenTnnn} for examples.

\begin{table}[h]
  \centering

  \begin{tabular}[h]{l|l|l|l|l}
    &\multicolumn{4}{c}{$j$ must satisfy these constraints}\\
    \cline{2-5}
    $n$ odd & $\tleft{j}$ & $\tright{j}$ & $\bleft{j}$ & $\bright{j}$ \\
    \hline
    \multirow{2}{*}{$\Tn{n}(\tleft{i})$}
    &$d_n(\tleft{j}) < d_n(\tleft{i})$ 
    & conflict 1(c)  
    & conflict 2(b)  
    & $d_n(\bright{j}) < d_n(\tleft{i}) $\\  
    &$ \Rightarrow j<i$  
    & $\Rightarrow i<j$  
    &$\Rightarrow j< \frac{n-3}{2} -i$  
    &$\Rightarrow  j> \frac{n-4}{2} -i$\\  
    \cline{2-5}

    \multirow{2}{*}{$\Tn{n}(\bleft{i})$}
    & conflict 2(b)  
    &$d_n(\tright{j}) < d_n(\bleft{i})$  
    & $d_n(\bleft{j}) < d_n(\bleft{i}) $  
    & conflict 1(c)\\  
    &$\Rightarrow j< \frac{n-3}{2} -i$  
    &$\Rightarrow j> \frac{n-4}{2} -i$  
    &$ \Rightarrow j<i$   
    & $\Rightarrow i<j$ \\  
    \cline{2-5}
    $\Tn{n}(\tright{i})$    & \multicolumn{4}{c}{symmetry with $\Tn{n}(\tleft{i})$}\\
    \cline{3-4}
    $\Tn{n}(\bright{i})$    & \multicolumn{4}{c}{symmetry with $\Tn{n}(\bleft{i})$}\\
    \cline{2-5}
    \multirow{2}{*}{$\Tn{n}(\varnothing)$}
    & $d_n(\tleft{j})<d_n(\tright{j})$ 
    & $d_n(\tright{j})<d_n(\tleft{j})$ 
    & $d_n(\bleft{j})<d_n(\bright{j})$ 
    & $d_n(\bright{j})<d_n(\bleft{j})$\\ 
    &$ \Rightarrow j< \frac{n-5}{4}$ 
    & $\Rightarrow j> \frac{n-5}{4}$ 
    &$\Rightarrow j< \frac{n-3}{4}$ 
    &$\Rightarrow j> \frac{n-3}{4}$\\\\ 
    $n$ even    & $\ldown{j}$&$\lup{j}$&$\rdown{i}$&$\rup{j}$\\
    \hline
    \multirow{2}{*}{$\Tn{n}(\ldown{i})$}
    &$d_n(\ldown{j}) < d_n(\ldown{i})$    
    & conflict 1(c) 
    & conflict 2(b)    
    & $d_n(\rup{j}) < d_n(\ldown{i}) $\\   
    &$ \Rightarrow j<i$   
    & $\Rightarrow i<j$   
    &$\Rightarrow j< \frac{n-2}{2}-i$   
    &$\Rightarrow j> \frac{n-3}{2}-i$\\   
    \cline{2-5}
    \multirow{2}{*}{$\Tn{n}(\lup{i})$}
    & conflict 1(c)    
    &$d_n(\lup{j}) < d_n(\lup{i})$    
    & $d_n(\rdown{j}) < d_n(\lup{i}) $   
    & conflict 2(b)\\   
    &$ \Rightarrow j<i$    
    & $\Rightarrow i<j$   
    &$\Rightarrow j< \frac{n-3}{2}-i$   
    &$\Rightarrow j> \frac{n-4}{2}-i$\\   
    \cline{2-5}
    $\Tn{n}(\rdown{i})$    & \multicolumn{4}{c}{symmetry with $\Tn{n}(\ldown{i})$}\\
    \cline{3-4}
    $\Tn{n}(\rup{i})$    & \multicolumn{4}{c}{symmetry with $\Tn{n}(\lup{i})$}\\
    \cline{2-5}
    \multirow{2}{*}{$\Tn{n}(\varnothing)$}
    & $d_n(\ldown{j})<d_n(\lup{j})$ 
    & $d_n(\lup{j})<d_n(\ldown{j})$ 
    & $d_n(\rdown{j})<d_n(\rup{j})$ 
    & $d_n(\rup{j})<d_n(\rdown{j})$\\ 
    &$ \Rightarrow j< \frac{n-3}{4}$  
    & $\Rightarrow j>\frac{n-3}{4}$  
    &$\Rightarrow j<\frac{n-3}{4}$ 
    &$\Rightarrow j>\frac{n-3}{4}$\\  
  \end{tabular}

  \caption{Each  $\Tn{n}(a)$ has exactly $\log_2(|\Tn{n}(a)|)$ independently
    flippable diagonals.  See Figures \ref{fig:oddTnnn}--\ref{fig:evenTnnn}
    for examples.}
  \label{tab:tnnnCases}
\end{table}

For odd $n$, the diagonals $a$ and $\bar a$ of the monomer on the
middle column are equal in length.  Thus is is not flipped at all in
$\Tn{n}(\varnothing)$, giving the number $|\Tn{n}(\varnothing)| =
2^{n-3}$, however, $a,\bar a \in A$, so $|A| = n-1$.

For even $n$, no such symmetry exists, so we have
$|\Tn{n}(\varnothing)| = 2^{n-2}$, and we must partition
$\Tn{n}(\varnothing)$ into halves to obtain $n$ classes of size
$2^{n-3}$.\qed

\begin{corollary}\label{cor:T(n)}
  The number of $n\times n$ tatami tilings is $2^{n-1}(3n-4)+2$.
\end{corollary}
\proof By Corollary 1 in \cite{EricksonRuskeySchurch2011}, we have
that $T(n,m) = 0$ when $m>n$.  Let $T(n) = \sum_{m\ge 0} T(n,m)$, so
that
    \begin{align*}
    T(n) = n2^{n-1} + \sum_{i = 1}^{\lfloor n/2 \rfloor}\left(
      (n-2i)2^{n-2i} + (n-2i+1)2^{n-2i+1}\right),
  \end{align*}
  and notice that the sum simplifies to
  \begin{align*}
    T(n) = n2^{n-1} + \sum_{i=1}^{n-1} i2^{i}.
  \end{align*}
  Now we use the fact that $2^k + 2^{k+1} + \cdots + 2^{n-1} = 2^n-1 -
  2^k+1 = 2^n-2^k$ to rearrange the sum.
  \begin{align*}
    T(n) =& n2^{n-1} + \sum_{i=1}^{n-1} i2^{i} \\
    =& n2^{n-1} + (n-1)2^n - \sum_{i=1}^{n-1} 2^i\\
    =& 2^{n-1}(3n-4)+2
  \end{align*}
\qed

By cross referencing the result of Corollary \ref{cor:T(n)} with the
On-Line Encyclopedia of Integer Sequences (\cite{Jovovic2005}), a
surprising correspondence with integer compositions of $n$ was
discovered (see sequence A027992). No proof is included in the
sequence reference so we have included one here which illuminates a
symmetry in this correspondence.

We restate the result of \tref{square} in terms of the distance of a
bidimer or vortex from the boundary of the grid.

\begin{lemma}\label{lem:distance} For $1\leq k\leq n-1$, the
  number of tatami tilings of the $n\times n$ grid with a vortex or
  bidimer a distance $(k+1)/2 $ from the boundary is $(n-k)2^{n-k}$.
\end{lemma}
  
\proof In the proof of \tref{square} it was shown that the number of
tilings of the $n\times n$ grid with $m$ monomers and a vortex is
$m2^m$, and with a bidimer it is $(m+1)2^{m+1}$. Let $d$ be the
minimum distance from the center of a feature to the boundary of the
grid. By \lref{distm} in the case of a bidimer $m=n-2d$ and in the
case of a vortex $m=2d+1$. Thus, in either case, the number of tilings
written in terms of $d$ is $(n-2d+1)2^{n-2d+1}$ or equivalently,
  $(n-k)2^{n-k}$ where $k=2d-1$.  \qed

\begin{proposition}
  \label{prop:sos}
  The number of $n\times n$ tatami tilings is equal to the sum of the
  squares of all parts in all compositions of $n$.
\end{proposition}

\begin{figure}[h]
  \centering
\includegraphics[width=\textwidth]{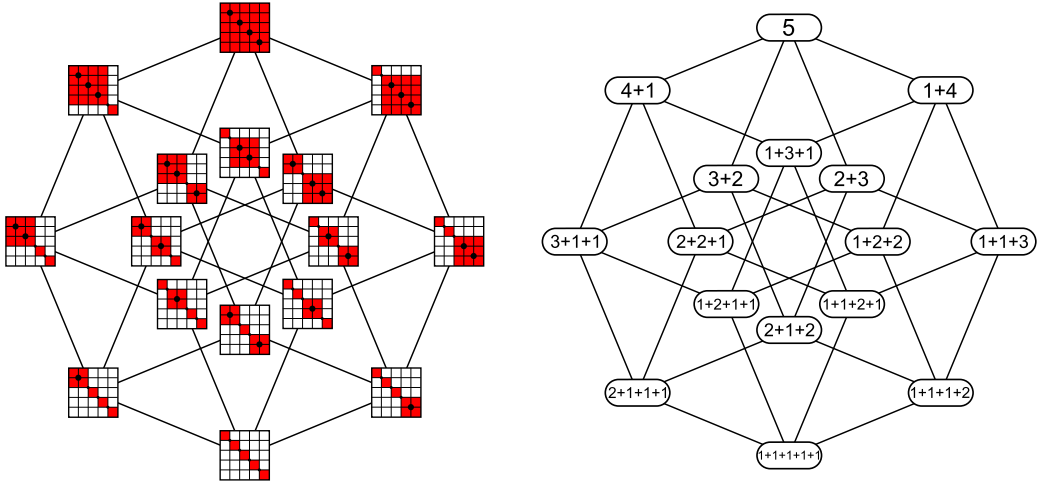}  
\caption{A pictorial representation of the squares of all parts in all
  compositions of $5$.}
  \label{fig:hass}
\end{figure}

\proof Note that the number of red boxes in the diagram in
\cite{Piesk2010}, see \fref{hass}, is the sum of the squares of all
parts in all compositions of $5$. We consider in general this visual
representation and show that for an integer composition of $n$ the
number of red boxes is equal to the number of tilings of the $n\times
n$ grid. We proceed by partitioning these red boxes into $n$ classes
and showing that the cardinality of each class is equal to the number
of tilings with a feature (or no features) a distance $k/2$ from the
edge of the grid, for some $k\in \{2,3,...,n\}$.

For each of the $n\times n$ arrays in \fref{hass}, view the top left
corner as having the $(x,y)$ coordinates $(0,0)$ and the bottom right
corner as having the coordinates $(n,-n)$. Now consider a grouping of
these red boxes into $n$ classes, where a box is in class $X_k$,for
$0\leq k\leq n-1$, whenever its center lies on one of the lines
$y=-x\pm k$.  The number of red boxes in class $X_0$ is the number of
red boxes along the line $y=-x$, which is $n2^{n-1}$ since there are
$2^{n-1}$ compositions of $n$ and $n$ red boxes along each central
diagonal. By \tref{tnnn} this is equal to the number of tilings of the
$n\times n$ grid with the maximum number of monomers, or equivalently,
no bidimers or vortices.

For $1\leq k \leq n-1$, if a summand of size $i\geq k$ appears in a
composition of $n$, then there are $2(i-k)$ red boxes along the lines
$y=\pm k$ and none if $i<k$.  In \cite{ChinnColyerFlashman1992} the
authors show that the number of summands of size $1\leq i <n$ over all
compositions of $n\geq 2$ is given by $(n+3-i)2^{(n-2-i)}$, and
trivially a summand of size $n$ occurs once. Thus the number of red
boxes in class $X_k$, for $1\leq k \leq n-1$,
is \begin{equation}\label{EQ:redsquaresum} 2(n-k)+\sum_{i=k+1}^{n-1}
  (n+3-i)(i-k)2^{(n-1-i)}.
\end{equation}

Expression (\ref{EQ:redsquaresum}) simplifies to $(n-k)2^{n-k}$
(details are omitted) which, by \lref{distance}, is the number of
tilings of the $n\times n$ grid with a feature a distance
$\frac{k+1}{2}$ from the boundary where $1\leq k\leq n-1$.  \qed

\section{Combinatorial generation algorithms}
\label{sec:combalg}


In this section we describe two combinatorial algorithms for
generating all elements of $\Tn{n}$, for even $n$.  Each algorithm is a natural
extension of our proof of \tref{tnnn}, and we use the notation and
ternary representation described there.

The first combinatorial generation algorithm arises most naturally
from our new proof of \tref{tnnn} (the proof in
\cite{EricksonRuskeySchurch2011} leads to another, recursive
algorithm).  For each $\Tn{n}(a)$, \tabref{tnnnCases} can be used to
determine the set of available diagonals. Several constant amortized
time (CAT) algorithms exist for generating all subsets, so our
algorithm is also CAT.

\subsection{Gray code}
\label{sec:gray}
The second algorithm uses a Gray code whose operation is the diagonal
flip, which, in our representation, is equivalent to incrementing or
decrementing a single ternary symbol.  We use a binary reflected Gray
code for each class of size $2^{n-3}$ tilings, and order the classes
so that the last tiling of one class is one flip away from first of
the next class.  If $n\equiv 2 \pmod 4$, this order begins with
\begin{align*}
\Tn{n}(\lup{0}),\Tn{n}(\rdown{\frac{n-2}{2}-1}),\Tn{n}(\lup{1}),
\Tn{n}(\rdown{\frac{n-2}{2}-2}), \ldots
,\Tn{n}(\lup{i}),\Tn{n}(\rdown{j})  
\end{align*}
where $i$ is the largest integer
such that $d_n(\lup{i}) > d_n(\ldown{i})$, and $j$ is the smallest
integer such that $d_n(\rdown{j})>d_n(\rup{j})$.  The next class is
$\Tn{n}(\varnothing)$, which is followed by
\begin{align*}
  \Tn{n}(\rup{0}),\Tn{n}(\ldown{\frac{n-2}{2}-1}), \Tn{n}(\rup{1}),
  \Tn{n}(\ldown{\frac{n-2}{2}-1}),\ldots \Tn{n}(\rup{j-1}),
  \Tn{n}(\ldown{i+1}).
\end{align*}
 For example, when $n=18$ (See \fref{tnnnRep}),
this order is
\begin{align*}
 \Tn{18}(\lup{0}),\Tn{18}(\rdown{7}),\Tn{18}(\lup{1}),\Tn{18}(\rdown{6}),
\Tn{18}(\lup{2}), \Tn{18}(\rdown{5}), \Tn{18}(\lup{3}),\Tn{18}(\rdown{4}),
\end{align*}
and then $\Tn{18}(\varnothing)$, continuing with
\begin{align*}
  \Tn{18}(\rup{0}),\Tn{18}(\ldown{7}), \Tn{18}(\rup{1}), \Tn{18}(\ldown{6}),
  \Tn{18}(\rup{2}), \Tn{18}(\ldown{5}), \Tn{18}(\rup{3}), \Tn{18}(\ldown{4}).
\end{align*}
 The case where $n\equiv 0 \pmod 4$ is
similar.

The binary reflected Gray code on $n$ bits begins with the $0$
bitstring and ends with $10^{n-1} $.  We assign $0$ to unflipped
monomers and $1$ to flipped monomers, since there is only one
available flip for each flippable monomer.  In this way, each
$\Tn{n}(a)$ can be generated by one of the CAT algorithms in
\cite{BitnerEhrlichReingold1976}.

To generate two consecutive classes, $\Tn{n}(a)$ and $\Tn{n}(b)$, they
must be separated by exactly one flip. Notice that if $b\neq
\varnothing$ and $a\neq \varnothing$, then $b$ is an available
diagonal in $\Tn{n}(a)$, since $a$ and $b$ are on different sides, in
different directions, and $d_n(b) = d_n(a)-1$.  In this case, we
associate $b$ with the $1$ in the last bitstring, $10^{n-3}$, of
$\Tn{n}(a)$, and unflipping $a$ allows us to begin generating
$\Tn{n}(b)$.  If $b = \varnothing$, then the order for the bitstring
of $\Tn{n}(a)$ is arbitrary, however, $\Tn{n}(\varnothing)$ must be
generated in reverse, ending with the trivial tiling.  If
$a=\varnothing$, then $\Tn{n}(\varnothing)$ ends with the trivial
tiling, so $b$ is flipped and $\Tn{n}(b)$ is generated.

The transitions from $\Tn{n}(a)$ to $\Tn{n}(b)$ are done in constant
time, so this Gray code generation algorithm for $\Tn{n}$ runs in
constant amortized time by making a sequence of calls to one of the
CAT algorithms for generating binary reflected Gray codes in
\cite{BitnerEhrlichReingold1976}.

Two problems remain open in this discussion.  Is there a Gray code for
$\Tn{n}$ whose first and last elements also differ by $1$ operation?
A positive answer would provide a Gray code for $\Tn{n}^{\text{all}}$.
 
Our present discussion is about $m=n$ and even $n$, but much of this
applies to tilings of odd dimension or with less than $n$ monomers.

\section{Further research}
\label{sec:conclusion}

This paper answers some questions central to the enumeration and
exhaustive generation of square monomer-dimer tatami tilings, and
similar findings are in preparation for the $r\times c$ grid, with
$r<c$ (see \cite{EricksonRuskeySchurch2011} for related conjectures).
The present research opens up the following paths:
\begin{itemize}
\item find a useful representation for tilings in $\Tn{n,m}$ which is
  ammenable to making diagonal flips;
\item exhaustively generate all of $\Tn{n,m}$ in constant amortized
  time;
\item find the corresponding Gray code of \sref{gray} for odd $n$;
\item find an aesthetically pleasing bijection from the $n\times n$
  tatami tilings to the red boxes of \fref{hass}, for all $n$.  Preferably one
  which uses elements of the proof of Proposition \ref{prop:sos}; and
\item explore $n\times n$ tilings with exactly $v$ vertical dimers.
\end{itemize}

This last point is in preparation, and might have been published here,
save for its length.  We will give a CAT algorithm for generating all
of the $n\times n$ tatami tilings with exactly $m$ monomers, $v$
vertical dimers and $h$ horizontal dimers, and we will explain the
origin of the generating polynomial in Conjecture 4 in
\cite{EricksonRuskeySchurch2011}.

\textbf{Acknowledgements:} The authors are grateful to Prof. Frank
Ruskey and Prof. Don Knuth for their interest and involvement in
tatami tilings.  We also thank the anonymous referees, and editors
Prof. Bill Smyth and Prof. Costas Iliopoulos.

\end{document}